\documentclass[11pt, a4paper, oneside, reqno]{amsart}
\allowdisplaybreaks
\usepackage[margin=1in]{geometry}
\usepackage{amsfonts,amsmath,amssymb,bbm,amsthm,stmaryrd}
\usepackage{amsthm}
\usepackage{hyperref}
\usepackage{mathtools}
\usepackage{enumerate}
\usepackage{xcolor}
\makeatletter
\renewcommand{\tocsection}[3]{%
	\indentlabel{\@ifnotempty{#2}{\bfseries\ignorespaces#1 #2.\,\,}}\bfseries#3}
\renewcommand{\tocsubsection}[3]{%
	\indentlabel{\@ifnotempty{#2}{\ignorespaces#1 #2\quad}}#3}
\renewcommand{\tocsubsubsection}[3]{%
	\quad\quad\quad\indentlabel{\@ifnotempty{#2}{\ignorespaces#1 #2\quad}}#3}

\newcommand\@dotsep{4.5}
\def\@tocline#1#2#3#4#5#6#7{\relax
	\ifnum #1>\c@tocdepth 
	\else
		\par \addpenalty\@secpenalty\addvspace{#2}%
		\begingroup \hyphenpenalty\@M
		\@ifempty{#4}{%
			\@tempdima\csname r@tocindent\number#1\endcsname\relax
		}{%
			\@tempdima#4\relax
		}%
		\parindent\z@ \leftskip#3\relax \advance\leftskip\@tempdima\relax
		\rightskip\@pnumwidth plus1em \parfillskip-\@pnumwidth
		#5\leavevmode\hskip-\@tempdima{#6}\nobreak
		\leaders\hbox{$\m@th\mkern \@dotsep mu\hbox{.}\mkern \@dotsep mu$}\hfill
		\nobreak
		\hbox to\@pnumwidth{\@tocpagenum{\ifnum#1=1\bfseries\fi#7}}\par
		\nobreak
		\endgroup
	\fi}
\AtBeginDocument{%
	\expandafter\renewcommand\csname r@tocindent0\endcsname{0pt}
}

\def\cM{\mathcal{M}}

\def\1{\mathbf{1}}

\newtheorem{theorem}{Theorem}[section]
\newtheorem{proposition}[theorem]{Proposition}
\newtheorem*{proposition*}{Proposition}
\newtheorem{corollary}[theorem]{Corollary}
\newtheorem*{theorem*}{Theorem}
\newtheorem{lemma}[theorem]{Lemma}

\theoremstyle{remark}
\newtheorem{remark}[theorem]{Remark}

\theoremstyle{definition}

\newcommand{\bnk}{\binom{[n]}{k}}

\title{Permutation Symmetries and characteristic polynomials}
\author{Nicolas Gilliers}
\begin{document}
\begin{abstract}
We compute certain multilinear functionals of permutation-invariant random matrices, following the work of Marcus and Geneson. As an application, we derive formulae for the average characteristic polynomial of the sum and the product of permutation-invariant random matrices.
\end{abstract}
\maketitle
\tableofcontents

\section{Introduction}
In these notes, we continue the work initiated by Marcus \cite{2016arXiv160706679M}
and Geneson \cite{2017arXiv170806640G}. In \cite{2016arXiv160706679M}, the author provides formulae for the mean characteristic polynomials of the sum and product of random matrices whose distributions are invariant under conjugation by a signed permutation matrix. Remarkably, the proposed formulae depend only on the characteristic polynomials of the original matrices. The associated convolutions are called \emph{finite free additive and multiplicative convolutions} and appeared ealier in the literature \cite{Walsh1922}. Invariance under signed permutation can be replaced by unitary invariance, and Marcus’s formulae remain valid; this phenomenon is related to the \emph{quadrature formula} enjoyed by unitary matrices and signed permutations. In recent years, much effort has focused on the asymptotic distributions of finite free convolutions of polynomials in connection with \emph{free probability} \cite{AGVP2023,ArizmendiPerales2018,Marcus2021Polynomial,Mirabelli2021,MSS2022FiniteFree}.
We prove formulae for the average characteristic polynomials of sums and products of permutation-invariant random matrices; see Corollaries \ref{cor:sum} and \ref{cor:mult}. As expected, the result depends not only on the spectral distribution of the initial matrices, but also on the distribution of the eigenspaces along the unique invariant subspace in the canonical representation of the permutation group as the Weyl group.

\subsection{Notations and basic formulae}
Let $n \geq 1 $ and $S_n \subset \cM_{n,n}(\mathbb{C})$ be the Weyl group. Given a matrix $M \in \mathcal{M}_{n,n}(\mathbb{C})$, we use $[M]_{X,Y}$ for the determinant of the matrix made of the columns in $X$ and rows in $Y$ of $M$. We set 
$$
[M]^{(k)} := \sum_{X \in \binom{[n]}{k}} [M]_{X,X}.
$$
With these notations, it holds that ${\rm det}(xI_n+M) = \sum_{k=0}^n [M]^{(k)}x^{n-k}$.
We denote by $\binom{[n]}{k}$ the set of all subsets of $[n]$ with cardinality $k$. Given two sets $X \in \binom{[n]}{k}$ and $Z \in \binom{[k]}{j}$ ($0\leq j\leq k)$, $X(Z)$ is the subsets of $[n]$ comprising all elements of $X$ in positions $Z$ in the increasing arrangement of the elements in $X$.
We recall the following basic formulae (which we use without mention)
\begin{align}
\tag{Cauchy-Binet}
&[M]_{X,Y} = \sum_{Z \in \binom{[n]}{k}} [M]_{X,Z}[M]_{Z,Y} \\
&[M]^{(n)} = \sum_{i=0}^n\sum_{X,Y \in \binom{[n]}{i}} (-1)^{\sum_{x\in X} x + \sum_{y\in Y} y} [M]_{X,Y}[M]_{X^c,Y^c} \nonumber
\end{align}
In the following we write ${\bf 1} = (1,\ldots,1)^T\in \mathbb{C}^n$. In the sequel, we set 
$$
|X|_1 := \sum_{x \in X}x.
$$
We denote by ${\rm adj}(M)$ the adjoint matrix of $M$, that is the transpose of the comatrix of $M$ :
$$
{\rm adj}(M)_{a,b}=(-1)^{a+b}[M]_{[n]-b, [n]-a}.
$$
\subsection{Contributions}
\begin{proposition} 
\label{prop:tocompute}
Let $0 \leq k\leq n$ two integers. Let $X,Y \in \binom{[n]}{k}$ two $k$-sets. Let $A,B,C \in \mathcal{M}_{n,n}(\mathbb{C})$ three matrices. Then
\begin{align*}
    &\sum_{P\in S_n}[BPAP^{-1}C]_{X,Y}= k!(n-k)![BC]_{X,Y}[A]^{(k)}\\
    &\hspace{2cm}+(k-1)!(n-k+1)!\Big(([B(I+{\bf 1}{\bf 1}^T)C]_{X,Y}-(k+1)[BC]_{XY}) \\ 
    &\hspace{3cm}\times ((n-k)[A]^{(k)}+\sum_{i=0}^k[A]^{(i)}{\rm Tr}(A^{j-i} {\bf 1}{\bf 1}^T))\Big)
\end{align*}
\end{proposition}
For any $A\in\mathcal{M}_{n,n}(\mathbb{C})$ and $1 \leq k\leq n$, set 
\begin{align}
\label{def:braces}
\{A\}_k:=((n-k)[A]^{(k)}-\sum_{i=0}^{k}(-1)^{k-i}[A]^{(i)} {\rm Tr}(A^{k-i}{\bf 1}{\bf 1}^T)
\end{align}
\begin{corollary} 
\label{cor:mult}
Let $A,B,C$ be as in proposition \ref{prop:tocompute}, and define $(c_k)_{n\geq k\geq 1}$ by:
\begin{align*}
\sum_{P\in S_n}{\rm det}(xI+BPAP^{-1}C) =  x^n+\sum_{k=1}^n c_k x^{n-k}
\end{align*}
Then, with $1 \leq k \leq n$, recalling \eqref{def:braces}
\begin{align*}
c_k&=\frac{k!(n-k)!}{n!}[CB]^{(k)}[A]^{(k)}+\frac{(k-1)!(n-k+1)!}{n!} \{A\}_k\{CB\}_k.
\end{align*}
\end{corollary}

\begin{corollary} 
\label{cor:sum}
Let $A,B$ be as in Proposition \ref{prop:tocompute} and define $(s_k)_{n\geq k \geq 1}$ by 
\begin{align}
&\sum_{P \in S_n} {\rm det}(xI+PAP^{-1}+B) = x^n+\sum_{k=1}^n d_kx^{n-k}
\end{align}
Then, with $1 \leq k \leq n$, recalling \eqref{def:braces} :
    \begin{align}
&d_k= \sum_{i=0}^k \frac{(n-i)!(n-k+i)!}{(n-k)!}[A]^{(i)}[B]^{(k-i)} + \sum_{i=1}^{k} \frac{(n-i+1)!(n-k+i-1)!}{(n-k)!} \{A\}_i\{B\}_{k-i}
\end{align}
\end{corollary}
\section{Proof of Proposition \ref{prop:tocompute}}
We start with a formula by J. Geneson. We need some notation. Given $X$ and $Y$ two subsets of $[n]$, let 
$$
s(X,Y)=\sum_{t\in X,t\not\in Y} r(t,X) +\sum_{t\in Y, t\not\in X}r(t,Y) = r(Y,X)+r(X,Y).
$$
where $r(t,X)$ is the number of elements of $X$ less than $t$. 
\begin{proposition}[\cite{2017arXiv170806640G}]
\label{prop:geneson}Let $S,T,U,V \in \binom{[n]}{k}$, then
\begin{equation*}
\sum_{P\in S_n} [P]_{S,T}[P^{-1}]_{U,V} = \begin{cases} |S\cap V|!(n-|S\cup V|)!(-1)^{s(S,V)+s(T,U)} & {\rm if}\, |S\cap V|=|U\cap T|\geq k-1 \\ 
0 & {\rm otherwise}
\end{cases}
\end{equation*}
\end{proposition}
The above formula relies on the simple observation that $[P]_{X,Y} \neq 0$ implies $Y=\sigma(X)$ and moreover, the matrix $(P)_{S,T}$( and $[P]_{U,V})$ can be brought to the following matrix, where $j = |S\backslash S\cap V|=|T\backslash U\cap T|$, 
$$
\begin{bmatrix}
I_{j} & 0 \\
0 &(P)_{U\cap S, V\cap T}
\end{bmatrix}
$$
by doing $s(S,V)$ swaps (lines or columns). We use the Cauchy-Binet formula to write 
\begin{align*}
\sum_{P}[BPAP^{-1}C]_{X,Y} &= \sum_{P\in S_n}\sum_{\substack{S,T,U,V \\ \in \binom{[n]}{k}}} [B]_{X,S}[P]_{S,T}[A]_{T,U}[P^{-1}]_{U,V}[C]_{V,Y} \\
&=\sum_{\substack{S,T,U,V \\ \in \binom{[n]}{k}}}[B]_{X,S}[A]_{T,U}[C]_{V,Y}\sum_{P \in S_n} [P]_{S,T}[P^{-1}]_{U,V}\\
&=\sum_{\substack{S,T,U,V \\ \in \binom{[n]}{k}\\ |S\cap V|=|U\cap T| \\ \geq k-1}}[B]_{X,S}[C]_{V,Y} [A]_{T,U} (-1)^{s(S,V)+s(T,U)}|S\cap V|!(n-|S\cup V|)! \\
&=k!(n-k)!\sum_{\substack{S,U \\ \in \binom{[n]}{k}}} [B]_{X,S}[C]_{S,Y} [A]_{U,U} \\ 
&\hspace{1cm}+(k-1)!(n-k+1)! \sum_{\substack{S,T,U,V \\ \in \binom{[n]}{k}\\ |S\cap V|=|U\cap T| \\ = k-1}}[B]_{X,S}[C]_{V,Y} [A]_{T,U} (-1)^{s(S,V)+s(T,U)}
\end{align*}
Let $S_1$ be the first sum on the right-hand side of the last expression, and set
\begin{align}
&S_3(B,C; X,Y)=  \sum_{\substack{S,V \\ \in \bnk \\ |S\cap V|=k-1}}[B]_{X,S}[C]_{V,Y} (-1)^{s(S,V)},\\ 
&S_2(A,k)= \sum_{\substack{T,U \\ \in \bnk \\|U\cap T|=k-1}}[A]_{T,U} (-1)^{s(T,U)}
\end{align}
so that
$$
\sum_{P \in S_n}[BPAP^{-1}C]_{X,Y}  = k!(n-k)!S_1+(k-1)!(n-k+1)!S_2(A,k)S_3(B,C;X,Y).
$$
\subsection{Computation of $S_1$} This sum has been computed in \cite{2016arXiv160706679M}
\begin{equation}
\label{eqn:sumssone}
S_1 = [BC]_{X,Y}[A]^{(k)}
\end{equation}
\subsection{Computation of $S_2$}
We write $S_2$ for brevity.
We rewrite the sum over $(U,T)$ as in $S_2$ as a the sum over $Z \in \binom{[n]}{k+1}$ and $(a,b)$ with $a\neq b \in [k+1]$:
\begin{align*}
S_2 &= \sum_{Z \in \binom{[n]}{k+1}}\sum_{a\neq b \in [k+1]} (-1)^{a+b-1} [A]_{Z([k+1]\backslash \{a\}),Z([k+1]\backslash \{b\})} \\
&= \sum_{Z \in \binom{[n]}{k+1}} \sum_{a,b \in [k+1]} (-1)^{a+b-1} [A]_{Z([k+1]\backslash \{a\}),Z([k+1]\backslash \{b\})}\\
&\hspace{2cm}+\sum_{Z\in\binom{[n]}{k+1}}\sum_{a\in[k+1]} [A]_{Z([k+1]\backslash \{a\}),Z([k+1]\backslash \{a\})} \\
&= \sum_{Z \in \binom{[n]}{k+1}} \sum_{a,b \in [k+1]} (-1)^{a+b-1} [A]_{Z([k+1]\backslash \{a\}),Z([k+1]\backslash \{b\})} \\ 
&\hspace{2cm}+(n-k)\sum_{W\in\binom{[n]}{k}}[A]_{W,W}
\end{align*}
To compute the first term in the last equation, we observe the following relation:
\begin{align}
\label{eqn:fund}
{\rm det}(M+{\bf 1}{\bf 1}^T)= \sum_{X,Y} (-1)^{|X|_1+|Y|_1}[M]_{X,Y}[{\bf 1}{\bf 1}^T]_{X^CY^C} = \sum_{a,b}[M]_{[n]\backslash a, [n]\backslash b} (-1)^{a+b}+[M]
\end{align}
We infer  
\begin{equation}
\label{eqn:plugging}
S_2 = \sum_{Z \in \binom{[n]}{k+1}}[A]_{Z,Z} - [A+{\bf 1}{\bf 1}^T]_{Z,Z} - (n-k)[A]^{(k)} = [A]^{(k+1)}-[A+{\bf 1}{\bf 1}^T]^{(k+1)}+(n-k)[A]^{(k)}
\end{equation}
Again, equation \eqref{eqn:fund} yields 
$$
[xI+A+{\bf 1}{\bf 1}^T]^{(n)} = [xI+A]^{(n)}+{\bf 1}^T{\rm adj}(xI+A){\bf 1}
$$
and we obtain $[A+{\bf 1}{\bf 1}^T]^{(k+1)}=[A]^{k+1}+(x^{n-1-k}){\bf 1}^T{\rm adj}(xI+A){\bf 1}$. We give a formula for the polynomial ${\rm adj}(xI+A)$ by using the fundamental relation :
$$
(xI+A){\rm adj}(xI+A)=[xI+A]^{(n)}I.
$$
Writing ${\rm adj}(xI+A)=\sum_{j=0}^{n-1}x^{n-1-j}C_j$, we obtain, whenever $n-1\geq j \geq 1$:
\begin{align*}
& [A]^{(j)}I= C_{j}+AC_{j-1} \\
&C_0 = I,
\end{align*}
and we infer :
$$
C_j = \sum_{i=0}^j(-1)^{j-i}[A]^{(i)}A^{j-i}.
$$
Finally,
\begin{align}
\label{eqn:rankone}
[A+{\bf 1}{\bf 1}^T]^{(k+1)}=[A]^{(k+1)}+\sum_{i=0}^k(-1)^{k-i}[A]^{(i)}{\rm Tr}(A^{k-i}{\bf 1}{\bf 1}^T).
\end{align}
Plugging \eqref{eqn:rankone} into \eqref{eqn:plugging} yields the following lemma :
\begin{lemma} 
\label{lem:stwoœ}
Let $A \in \mathcal{M}_{n}(\mathbb{C})$ and $1\leq k \leq n$, then
$$
S_2(A,k) = (n-k)[A]^{(k)}+\sum_{i=0}^k(-1)^{i-1}[A]^{(k-i)}{\rm Tr}(A^{i} {\bf 1}{\bf 1}^T)
$$
\end{lemma}
\begin{remark}
If $A$ is diagonal, then $S_2$ is equal to $0$. In fact, if $A = {\rm diag}(a_1,\ldots,a_n)$, then 
\begin{equation}
\label{eqn:newtons}
\sum_{i=0}^k(-1)^{i-1}[A]^{(k-i)}{\rm Tr}(A^{i} {\bf 1}{\bf 1}^T) = -ne_k(a_1,\ldots,a_n) + \sum_{i=1}^k(-1)^{i-1}e_{k-i}(a_1,\ldots,a_n)p_i(a_1,\ldots,a_n)
\end{equation}
where the $e_i(a_1,\ldots,a_n)$ and $p_i(a_1,\ldots,a_n)$, $n\geq i \geq 1=$ are the elementary symmetric polynomials and the Newton's sums of the $a's$. By Newton's identities, the right-hand side of \eqref{eqn:newtons} is trivial. Hence the equation $\{A\}_k = 0$, $1 \leq k \leq n$ is a matricial equivalent of the Newton's identities.
\end{remark}
\subsection{Computation of $S_3$}
As for $S_2$, we will rewrite the sum over $S$ and $V$ as a double sum over $W \in \binom{[n]}{k-1}$ and $a,b \in [n]\backslash W$ : 
\begin{align*}
S_3 &= \sum_{W \in \binom{[n]}{k-1}}\sum_{a\neq b \in [n]\backslash }[B]_{X,W \cup a)}[C]_{W\cup b,Y}(-1)^{s(W\cup a, W\cup b)} \\
&= \sum_{W \in \binom{[n]}{k-1}}\sum_{a\neq b \in [n]\backslash W }[(B)_{X,a} | (B)_{X,W}]\begin{bmatrix}  (C)_{b,Y}\\ (C)_{W,Y} \\\end{bmatrix} \\
&= \sum_{W \in \binom{[n]}{k-1}}\sum_{a, b \in [n] }[(B)_{X,a} | (B)_{X,W}]\begin{bmatrix}  (C)_{b,Y}\\ (C)_{W,Y} \\\end{bmatrix} \\
&\hspace{2cm}-\sum_{W \in \binom{[n]}{k-1}}\sum_{a\in [n] }[(B)_{X,a} | (B)_{X,W}]\begin{bmatrix}  (C)_{a,Y}\\ (C)_{W,Y} \\\end{bmatrix} \\
&= \sum_{W \in \binom{[n]}{k-1}}[(B1)_{X} | (B)_{X,W}]\begin{bmatrix}  ({\bf 1}^TC)_{Y}\\ (C)_{W,Y} \\\end{bmatrix} \\
&\hspace{2cm}-\sum_{W \in \binom{[n]}{k-1}}\sum_{a\in [n] }[(B)_{X,a} | (B)_{X,W}]\begin{bmatrix}  (C)_{a,Y}\\ (C)_{W,Y} \end{bmatrix} \\
&=\sum_{W\in \binom{\{0\}\cup [n]}{k}}[B1 | B]_{X,W}\begin{bmatrix}  {\bf 1}^TC\\ C \\\end{bmatrix}_{W,Y} -\sum_{W\in \binom{[n]}{k-1}}[B]_{X,W}[C]_{W,Y}\\
&\hspace{2cm}-\sum_{W \in \binom{[n]}{k-1}}\sum_{a\in [n] }[(B)_{X,a} | (B)_{X,W}]\begin{bmatrix}  (C)_{a,Y}\\ (C)_{W,Y} \end{bmatrix} \\
\end{align*}
On the first line of the last equation, we recognise the Cauchy-Binet formula for $[B({\bf 1}|I)\begin{pmatrix}{\bf 1}^T \\ I\end{pmatrix}C]_{X,Y}=[B(I+J)C]_{X,Y}$. On the second line, we move the first column of the first matrix and the first line of the second matrix to obtain 
$$
\sum_{W \in \binom{[n]}{k-1}}\sum_{a\in [n] }[B]_{X,W\cup\{a\}}\begin{bmatrix} C_{W\cup \{a\},Y} \end{bmatrix} = k [BC]_{X,Y}
$$
by the Cauchy-Binet formula. We obtain 
\begin{lemma} 
\label{lem:sthree}
Let $B,C \in  \mathcal{M}_n(\mathbb{C})$ and $X,Y \in \bnk$, then
\begin{equation*}
S_3(B,C; X,Y)=[B(I+{\bf 1}{\bf 1}^T)C]_{X,Y} - (1+k)[BC]_{X,Y}.
\end{equation*}
\end{lemma}
Equation \eqref{eqn:sumssone}, lemmas \ref{lem:stwoœ} and \ref{lem:sthree} imply Proposition \ref{prop:tocompute}.
\section{Proof of Corollary \ref{cor:mult}}
We start with the following identity:
\begin{equation*}
{\rm det}(xI+BPAP^{-1}C)=\sum_{X\in\binom{[n]}{k}}\sum_{P}[BPAP^{-1}C]_{X,X}
\end{equation*}
and insert the formula stated in Proposition \ref{prop:tocompute}. We evaluate first
\begin{align*}
[B(I+{\bf 1}{\bf 1}^T)C]_{X,X}=\sum_{S,T \in \bnk} [B]_{X,S}[I+{\bf 1}{\bf 1}^T]_{S,T}[C]_{T,Y}
\end{align*}
Observe that, if $|S\cap T|=k-1$ then $[I+{\bf 1}{\bf 1}^T]_{S,T} = (-1)^{s(S,T)}[{\rm diag}(0,I_{k-1})+{\bf 1}{\bf 1}^T]=(-1)^{s(S,T)}$. If $|S\cap T| \leq k-2$, then $[I+{\bf 1}{\bf 1}^T]_{S,T}$ contains two columns equal to ${\bf 1}$ and therefore $[I+{\bf 1}{\bf 1}^T]_{S,T}=0$. If $S=T$, $[I+{\bf 1}{\bf 1}^T]_{S,T}=k+1$. Thus 
$$
[B(I+{\bf 1}{\bf 1}^T)C]_{X,X}=(k+1)[BC]_{X,X}+\sum_{\substack{S,T\in\bnk \\ |S\cap T|=k-1}} (-1)^{s(S,T)} [B]_{X,S}[C]_{T,X}.
$$
Summing over $X \in \bnk$ we get 
\begin{align*}
\sum_{X \in \bnk}[B(I+{\bf 1}{\bf 1}^T)C]_{X,X}&=(k+1)[BC]^{(k)}+\sum_{\substack{S,T \in \bnk \\ |S\cap T|=k-1}} (-1)^{s(S,T)}[CB]_{T,S} \\
&= (k+1)[BC]^{(k)}+(n-k)[CB]^{(k)}+\sum_{i=0}^k(-1)^{k-i}[CB]^{(i)}{\rm Tr}((CB)^{k-i} {\bf 1}{\bf 1}^T) \\
&=(n+1)[BC]^{(k)}+\sum_{i=0}^k(-1)^{k-i}[BC]^{(i)}{\rm Tr}((CB)^{k-i} {\bf 1}{\bf 1}^T)
\end{align*}
The result follows.

\section{Proof of Corollary \ref{cor:sum}}
We use the Geneson's formula (Proposition \ref{prop:geneson}) to write 
\begin{align*}
&\sum_{P \in S_n} [PAP^{-1}+B]^{(k)} =\sum_{P \in S_n} \sum_{X} \sum_{i=0}^k \sum_{S,T \in \binom{[k]}{i}}(-1)^{|S|_1+|T|_1}[PAP^{-1}]_{X(T),X(S)}[B]_{X(S^c),X(T^c)} \\
& =\sum_{X\in \bnk} \sum_{i=0}^k  i!(n-i)![A]^{(i)} [B]_{X(S^c),X(T^c)}\\ 
&+ (i-1)!(n-i+1)!\sum_{X} \sum_{i=1}^k \sum_{S,T \in \binom{[k]}{i}}(-1)^{s(X(S),X(T))}(-1)^{|S\Delta T|_1}\delta_{|S\cap T|=i-1}\\
&\hspace{6cm}((n-i)[A]^{(i)}+\sum_{l=0}^i(-1)^{i-l}[A]^{(l)} {\rm Tr}(A^{i-l}11^T) [B]_{X(S^c),X(T^c)}
\end{align*}
At this stage, we have to find a formula for the sign 
$$
(-1)^{s(X(S),X(T)) + |S\Delta T|_1}
$$
Write $S = W^c \cup \{s\}$ and $T = W^c \cup \{t\}$. By definition, $s(X(S),X(T)) = |\{ u \in W^c : u < s\}| + \{ u \in W^c : u < t\}|$, which we rewrite as : 
$$
s(X(S),X(T)) = s-1 - |\{u \in W : u < s \} + t-1 - |\{u \in W : u < t \} = |S\Delta T|_1 + r(s,W)+r(t,W). 
$$
We obtain:
\begin{align*}
&= \sum_{i=0}^k \frac{(n-i)!(n-k+i)!}{(n-k)!}[A]^{(i)}[B]^{(k-i)}\\
&\hspace{0.5cm}+ (i-1)!(n-i+1)!\sum_{X} \sum_{i=1}^k ((n-i)[A]^{(i)}+\sum_{l=0}^i[A]^{(l)} {\rm Tr}(A^{i-l}11^T)) \\ 
&\hspace{8cm}\sum_{W \in \binom{k}{k-i+1}}\sum_{s\neq t \in W} (-1)^{r(s,W)+r(t,W)}[B]_{X(W-s),X(W -t)} \\
&= \sum_{i=0}^k \frac{(n-i)!(n-k+i)!}{(n-k)!}[A]^{(i)}[B]^{(k-i)}\\
&+ \sum_{X \in \bnk} \sum_{i=1}^k (i-1)!(n-i+1)!((n-i)[A]^{(i)}+\sum_{l=0}^i(-1)^{i-l}(-1)^{l}[A]^{(l)} {\rm Tr}(A^{i-l}11^T))S_2((B)_{X,X},k-i+1)
\end{align*}
We use the relations 
$$
S_{2}((B)_{X,X}, k-i+1)= [(B)_{X,X}]^{(k-i+1)}-[(B+1^T1)_{X,X}]^{(k-i+1)}+i[(B)_{X,X}]^{(k-i)}
$$
Now if we sum over $X \in \bnk$ the last expression, we obtain : 
\begin{align*}
\sum_{X\in \bnk} S_{2}((B)_{X,X}, k-i) &= \binom{n-k+i-1}{i-1}([B]^{(k-i+1)}-[(B+1^T1)_{X,X}]^{(k-i+1)}) \\
&\hspace{6cm}+ i\binom{n-k+i}{i} [(B)_{X,X}]^{(k-i)} \\
&= \binom{n-k+i-1}{i-1}([B]^{(k-i+1)}-[(B+{\bf 1}{\bf 1}^T)_{X,X}]^{(k-i+1)} \\
&\hspace{6cm}+ ({n-k+i})[(B)_{X,X}]^{(k-i)})
\end{align*}
Equation \eqref{eqn:rankone} yields the result.
\bibliographystyle{alpha}
\bibliography{main}

\end{document}